\documentclass[11pt]{amsart}

\usepackage{graphicx}
\usepackage{amscd}
\usepackage{amsmath}
\usepackage{amsfonts}
\usepackage{amssymb}
\usepackage{amsxtra}
\usepackage{xspace}
\usepackage{array}
\usepackage{cite}
\usepackage{setspace}
\theoremstyle{plain}
\newtheorem{theorem}{Theorem}

\newtheorem{corollary}{Corollary}
\newtheorem{proposition}{Proposition}

\theoremstyle{definition}

\numberwithin{equation}{section}

\newcommand{\be}{\begin{enumerate}}
\newcommand{\ee}{\end{enumerate}}
\newcommand{\beq}{\begin{equation}}
\newcommand{\eeq}{\end{equation}}
\newcommand{\bprop}{\begin{proposition}}
\newcommand{\eprop}{\end{proposition}}
\newcommand{\complex}{\mathbb{C}}

\newcommand{\integers}{\mathbb{Z}}

\DeclareMathOperator{\len}{\ell}

\DeclareMathOperator{\ct}{ct} 
\DeclareMathOperator{\Max}{Max} 
\newcommand{\s}{\scriptscriptstyle}
\newcommand{\Thetaram}{\Theta_R}
\newcommand{\Thetaone}{\Theta_1}
\newcommand{\Thetatwo}{\Theta_2}
\newcommand{\Thetatwobar}{\widehat{\Theta}_2}
\newcommand{\Thetafour}{\Theta_4}
\newcommand{\alm}[1][t,q]{a^{\s \lambda}_{\s \mu}(#1)}

\newcommand{\kost}[1][\mu]{K_{\lambda #1}(t)}
\newcommand{\hl}[1][\lambda]{P_{#1}(t)}
\newcommand{\sch}[1][\lambda]{\chi_{{}_{#1}}}

\newcommand{\Ep}{\mathcal{E}}

\newcommand{\kma}{\mathfrak{g}}
\newcommand{\csa}{\mathfrak{h}}

\newcommand{\dtl}{\tilde{\varDelta}}

\newcommand{\dtli}{\tilde{\varDelta}^{im}}
\newcommand{\aoo}{\ensuremath{a^0_0(t,q)}\xspace}

\newcommand{\alolo}[1][t,q]{\ensuremath{a^{\Lambda_0}_{\Lambda_0}(#1)}\xspace}
\newcommand{\kmaf}{\overset{\circ}{\mathfrak{g}}}

\newcommand{\chke}{\hat{\boldsymbol{\mu}}}

\newcommand{\tsf}{$t$-string function\xspace}
\newcommand{\tsfs}{$t$-string functions\xspace}

\newcommand{\Ept}{\Ep_t}
\newcommand{\ro}{\Delta}
\newcommand{\rp}{\Delta_+}
\newcommand{\rmin}{\Delta_-}
\newcommand{\rre}{\Delta^{re}}
\newcommand{\rrep}{\Delta^{re}_+}
\newcommand{\rim}{\Delta^{im}}
\newcommand{\rimp}{\Delta^{im}_{+}}

\newcommand{\skma}{symmetrizable Kac-Moody algebra\xspace}

\newcommand{\afsl}{\ensuremath{A_1^{(1)}}}
\newcommand{\ea}[1][\alpha]{e^{- #1}}
\newcommand{\hyp}[2][\infty]{(#2)_{#1}}
\newcommand{\nh}[1][]{[#1]_\infty}
\newcommand{\ps}{\mathbb{F}}
\newcommand{\fl}[2][]{F_{#1, #2}(t,v)}
\newcommand{\bil}[1][]{{}_{#1}\psi_{#1}}

\begin{document}

\title[]{Affine Hall-Littlewood functions for $\afsl$ and some 
 constant term identities of Cherednik-Macdonald-Mehta type}
\author{Sankaran Viswanath}
\address{Department of Mathematics\\
Indian Institute of Science\\
Bangalore 560012, India.}
\email{svis@math.iisc.ernet.in}
\thanks{\scriptsize Supported by UGC-SAP/DSA-IV}

\begin{abstract}
We study $t$-analogs of string functions for integrable highest weight 
representations of the affine Kac-Moody algebra
$\afsl$. We obtain closed form formulas for certain 
$t$-string functions of levels 2 and 4. As corollaries, we obtain
explicit identities for the corresponding affine Hall-Littlewood
functions, as well as higher-level generalizations of Cherednik's Macdonald and
Macdonald-Mehta constant term identities.
\end{abstract}

\maketitle

\singlespace
\section{Introduction}

Let $\kma$ be an untwisted affine Kac-Moody algebra. Let $\delta$
denote its null root, $P^+$ denote the set of dominant integral
weights and let $t$ be an indeterminate. To each $\lambda \in P^+$,
one can associate a Hall-Littlewood function $\hl$; these are affine
generalizations of the classical Hall-Littlewood polynomials
\cite{ramnelsen} and interpolate between the Weyl-Kac characters
$\sch$ and the Weyl group orbit sums $m_\lambda$
\cite{groj,svis-kosfoulkes}.  The affine Kostka-Foulkes polynomials
$\kost$ are defined to be the coefficients that appear in the relation
$\sch = \sum_{\mu \in P^+} K_{\lambda\mu}(t) P_{\mu}(t)$.  Affine
Hall-Littlewood functions have been studied in various contexts by
several authors \cite{groj,fgt,kap,cher-new,svis-kosfoulkes}. 
Affine Kostka-Foulkes polynomials are also of great interest,
especially on account of their positivity properties \cite{fgt}.

One has a well developed theory of integrable highest weight
representations of affine Kac-Moody algebras.  Given $\lambda \in
P^+$, let $L(\lambda)$ denote the corresponding highest weight
representation.  The generating functions for weight multiplicities
along $\delta$-strings in $L(\lambda)$ gives (upto a multiplicative
factor) the Kac-Moody string
functions \cite[(12.7.7)]{kac}. Since affine Hall-Littlewood functions
$\hl$ are $t$-deformations of the $m_\lambda$, the coefficients
$\kost$ above may be viewed as $t$-analogs of weight
multiplicities. The notion of a string function thus admits a natural
$t$-analog; these so called \tsfs \cite{groj,svis-kosfoulkes} are
defined to be generating functions for the $\kost$ along
$\delta$-strings.

It is known (see \cite{cher-new} for an extensive account) that affine
Hall-Littlewood functions are closely related to the theory of the
{\em Double Affine Hecke algebra} (DAHA). In fact, our notion is a
special case (corresponding to dominant affine weights) of the more
general notion of affine {\em Hall functions of level} $l>0$
introduced in \cite{cher-new}. The level $1$ case has been studied for
all affine root systems \cite{dmmc,stokman,cher-new}; for simply laced
$\kma$ and $\lambda$ of level 1, there are closed form 
formulas for $\hl$ and the corresponding \tsf. In fact, the explicit formula for the level
1 \tsf is equivalent to  
Cherednik's Macdonald-Mehta constant term identity \cite{dmmc, svis-kosfoulkes}.
However, there is almost nothing known about higher level affine
Hall-Littlewood functions.

In this article, we consider the rank 1 affine Kac-Moody algebra
$\afsl$.  We study higher level \tsfs
functions with a view toward obtaining explicit formulas. We exhibit
such formulas for certain level 2 and 4 \tsfs; in turn these determine
the corresponding affine Hall-Littlewood functions.  Interpreted in
terms of the {\em constant term} functional \cite{dmmc}, our results
provide new constant term identities generalizing those of Macdonald
and Macdonald-Mehta (in the $\afsl$ case).

Our method is based on studying {\em principal specializations} of
affine Hall-Littlewood functions. When $\kma=\afsl$, these can be
related to certain bilateral $q$-hypergeometric series. For our
special dominant weights, the resulting series can be summed in closed
form using various specializations of Bailey's classical $\bil[6]$
summation.

Our methods are very specific to rank one and we do not know if they
can be extended to general untwisted affines, or even to
$A_n^{(1)}$. Nevertheless, our results provide the first explicitly
computable examples of higher level affine Hall-Littlewood functions
and we hope they will be of value in studying the role of the DAHA in this theory.

The article is arranged as follows. In \S \ref{prelimsec}, we recall
the main aspects of the theory of Hall-Littlewood functions
associated to affine (and more generally symmetrizable Kac-Moody) Lie
algebras. In \S \ref{prspsec}, we work with $\afsl$, and compute
principal specializations of Hall-Littlewood functions using
classical summation formulas for our special cases. Section
\ref{mainthms} contains our main results; specifically theorems
\ref{al0l1thm}-\ref{alev4thm} compute closed form
expressions for the \tsfs and corollaries
\ref{ct-lev2-1}-\ref{ct-lev4} give formulas  
for the affine Hall-Littlewood functions and the new constant term identities.
The short appendix collects
together classical summation theorems for $q$-hypergeometric series that
are used in this article.

\section{Hall-Littlewood functions}\label{prelimsec}
\subsection{}
We first recall relevant facts and notation concerning Hall-Littlewood
functions associated to symmetrizable Kac-Moody algebras
\cite{svis-kosfoulkes, groj, fgt}.  Let $\kma$ be a \skma with Cartan subalgebra
$\csa$. Let $\alpha_i, i=1 \cdots n$ be the simple roots of $\kma$ and
$P$, $Q$, $P^+$, $Q^+$ be the weight lattice, the root lattice and the
sets of dominant weights and non-negative integer linear combinations
of simple roots respectively. We will denote the Weyl group of $\kma$
by $W$ and let $(\cdot,\cdot)$ be a nondegenerate, $W$-invariant
symmetric bilinear form on $\csa^*$.  Let $\len(\cdot)$ be the length
function on $W$ and for $\lambda \in \csa^*$, let $W_\lambda \subset
W$ denote the stabilizer of $\lambda$.  Let $\ro$ (resp. $\rp$,
$\rmin$) be the set of roots (resp. positive, negative roots) of
$\kma$. We will also let $\rre$ (resp. $\rre_\pm$) denote the set of
real roots of $\kma$ (resp. positive/negative real roots), and
similarly $\rim$ (resp. $\rim_\pm$), the corresponding subsets of
imaginary roots.  Given $\lambda, \mu \in \csa^*$, we say $\lambda
\geq \mu$ if $\lambda - \mu \in Q^+$. Given $\lambda \in \csa^*$,
define $D(\lambda) :=\{ \gamma \in \csa^*| \gamma \leq \lambda\}$. Let
$\Ept$ be the set of all series of the form
\begin{equation}\label{ser}
\sum_{\lambda \in \csa^*} c_\lambda(t) \, e^{\lambda}
\end{equation}
where each $c_\lambda(t) \in \complex[[t]]$ and $c_\lambda =0$ outside the
union of a finite number of sets of the form $D(\mu), \, \mu \in \csa^*$.
 For each $\alpha \in \ro$, let   $\mathrm{mult} (\alpha)$ denote the root
multiplicity of $\alpha$. 
Let $\rho$ be a Weyl vector of $\kma$ defined by
 $(\rho, \alpha^\vee_i) =1 \; \forall i=1\cdots n$, where as usual
 $\alpha^\vee_i = \frac{2\alpha_i}{(\alpha_i,\alpha_i)}$. Finally,
 given $\lambda \in P^+$, let $L(\lambda)$ be the integrable
 $\kma$-module with highest weight $\lambda$.

 For $\lambda \in P^+$, the {\em Hall-Littlewood function} $\hl$
 is defined as \beq\label{objint} \hl := \frac{1}{W_\lambda(t)}
 \frac{\sum_{w \in W} (-1)^{\len(w)} w\left(e^{\lambda+\rho}
   \prod_{\alpha \in \rp} \left( 1-t e^{-\alpha}
   \right)^{\mathrm{mult} (\alpha)}\right)}{e^{\rho} \prod_{\alpha \in \rp} (1 -
   e^{-\alpha})^{\mathrm{mult} (\alpha)}} \eeq where ${W_\lambda(t)} =
 \sum_{\sigma \in W_\lambda} t^{\len(\sigma)}$ is the Poincar\'{e}
 series of $W_\lambda$. Proposition 1 of \cite{svis-kosfoulkes} implies
 that $\hl$ is a well defined element of $\Ept$.  If $\sch$ denotes
 the formal character of the irreducible highest weight representation
 $L(\lambda)$ of $\kma$, we can write 
 \beq\label{chikpl} \sch = \sum_{\mu \in P^+, \, \mu \leq \lambda} K_{\lambda\mu}(t) P_{\mu}(t).
 \eeq 
 By corollary 1 of \cite{svis-kosfoulkes}, the {\em Kostka-Foulkes
   polynomials} $ K_{\lambda\mu}(t)$ lie in $\integers[t]$.
One also has  $P_\lambda(0) = \sch$ and $P_\lambda (1) = m_{\s \lambda} :=\sum_{\mu
  \in W\lambda} e^\mu$. As a consequence, $K_{\lambda \mu}(0) =
\delta_{\lambda, \mu}$ and  $K_{\lambda \mu}(1) = \dim(L(\lambda)_\mu)$.

Define the special element $\dtl \in \Ept$ by
$$\dtl := \frac{\prod_{\alpha \in \rp}
    (1-e^{-\alpha})^{\mathrm{mult} (\alpha)}}{\prod_{\alpha \in \rp}
      (1-t e^{-\alpha})^{\mathrm{mult} (\alpha)}} =  \prod_{\alpha \in \rp}
    \left[(1-e^{-\alpha})(1+t e^{-\alpha} + t^2 e^{-2\alpha} +
      \cdots)\right]^{\mathrm{mult} (\alpha)}$$
Letting $S(w):=\{\alpha \in \rp: w^{-1} \alpha \in \rmin\}$, we have
the following useful relation 
(equation (4.2) of \cite{svis-kosfoulkes}).
\begin{equation} 
e^{-\lambda} \dtl \hl = \frac{1}{W_\lambda(t)} \sum_{w \in W}
e^{w\lambda - \lambda} \prod_{\alpha \in S(w)}
\frac{t-e^{-\alpha}}{1-te^{-\alpha}} \label{dlpl}
\end{equation}

\subsection{}
Specializing to the case that $\kma$ is an untwisted affine Kac-Moody
algebra, we let $\kmaf$ denote the underlying finite dimensional
simple Lie algebra of rank $l$, say. If $\delta$ is the null root,
then $\rimp =\{k\delta: k \geq 1\}$ and each imaginary root has
multiplicity $l$.

Suppose $\lambda \in P^+$, let $\Max(\lambda):=\{\mu \in P^+: \mu \leq
\lambda; \, \mu + \delta \nleq \lambda\}$. For each $\mu \in
\Max(\lambda)$, the generating function \beq \alm := \sum_{k \geq 0}
K_{\lambda, \mu-k\delta}(t) \; q^k \eeq is termed a {\em $t$-string
  function}.  The {\em constant term} map $\ct(\cdot)$ \cite{igmgod}
is defined on elements $f = \sum_{\lambda} f_{\lambda} e^\lambda$ of
$\Ept$ by $$\ct(f):=\sum_{k \in \integers} f_{k\delta} \:
e^{k\delta}.$$ If $g \in \Ept$ is such that $g=\ct(g)$, then $\ct(fg) =
g \ct(f)$ for all $f \in \Ept$ \cite[(3.1)]{igmgod}.

We will let $q:=e^{-\delta}$ for the rest of the paper. The key
relationship between the notions of the preceding paragraph is the following
 \cite[(5.8)]{svis-kosfoulkes}:
\beq\label{almeq}
\alm = \ct(e^{-\mu}\, \dtl\, \sch)
\eeq

\section{$\afsl$ and bilateral $q$-hypergeometric series} \label{prspsec}
For the rest of the article we will restrict ourselves to the simplest
affine Kac-Moody algebra $\afsl$. The goal of this section is to
compute the principal specializations of certain Hall-Littlewood
functions associated to $\afsl$.  Let $\alpha_0$, $\alpha_1$ denote
the simple roots of $\afsl$; the null root $\delta =\alpha_0 +
\alpha_1$. \onehalfspace The real roots $\{\alpha_i + k \delta: k \in \integers,
i=0,1\}$ and imaginary roots $\{k\delta: k \in \integers \backslash
\{0\}\}$ all have multiplicity 1.  We will assume that the form $(,)$
is normalized such that $(\alpha_i, \alpha_i) = 2$ for $i=0, 1$.
The Weyl group $W$ is the infinite dihedral group
generated by the simple reflections $r_0, r_1$; we write $W:=\{w_j: j
\in \integers\}$ where $w_j:=(r_0 r_1)^{\frac{j}{2}}$ for $j$ even and
$w_j:=r_0 (r_1 r_0)^{\frac{j-1}{2}}$ for $j$ odd. It is easily seen
that $S(w_0) = \emptyset$ and $S(w_j)$ equals $\{\alpha_0,
\alpha_0+\delta,\cdots,\alpha_0+(j-1)\delta\}$ for $j>0$ and
$\{\alpha_1, \alpha_1+\delta,\cdots,\alpha_1+(|j|-1)\delta\}$ for
$j<0$. For all $j \in \integers$ one has $\rho - w_j \rho =
\sum_{\gamma \in S(w_j)} \gamma = j \alpha_0 + {j \choose 2} \delta$.

Let $\pi(w_j):=\prod_{\alpha \in S(w_j)} \frac{t-\ea}{1-t\ea}$.
One observes \cite[(4.2)]{igmgod} that $$\pi(w_j) = t^{|j|}
\prod_{\alpha \in S(w_j)} \frac{1-t^{\s -1} \ea}{1-t\ea} = t^j \, \frac{(t^{\s -1} \ea[\alpha_0] ; \ea[\delta])_j}{(t \ea[\alpha_0] ; \ea[\delta])_j} $$
\onehalfspace
for all $j \in \integers$. We have used the usual $q$-hypergeometric
notations:  $\hyp{a; \, w} := \prod_{i=0}^\infty (1-a\, w^i)$ and  
 $(a;w)_j := (a; w)_\infty/ (aw^j; w)_\infty$ for all $j \in \integers$. We will also use
the shorthand $(a_1, a_2, \cdots, a_k; w)_j:= (a_1;w)_j \, (a_2;w)_j \, \cdots (a_k;w)_j$.

\singlespace
Next, for $\gamma = c_0 \alpha_0 + c_1 \alpha_1 \in Q$, let 
$\mathrm{ht} (\gamma) := c_0 + c_1 = (\gamma,\rho)$.  Given $f
:=\sum_{\beta \in Q^+} f_\beta \,\ea[\beta]$ with $f_\beta \in
\complex[[t]]$, the principal specialization $$\ps(f) := \sum_{\beta
  \in Q^+} f_\beta \, v^{\mathrm{ht} (\beta)} = \sum_{\beta \in Q^+}
f_\beta \, v^{ (\beta,\rho)} \in \complex[[t,v]].$$ For instance,
\begin{equation}
\ps(\dtl) = \frac{\hyp{v, v, v^2; v^2}}{\hyp{tv, tv, tv^2; v^2}} \text{ and
} \ps(\pi(w_j)) = t^j \frac{\hyp[j]{t^{\s -1} v ; v^2}}{\hyp[j]{t v ;
      v^2}}. \label{psdtl}
\end{equation}

\noindent
Since $\ps(\ea[\alpha_0]) = \ps(\ea[\alpha_1])$, we have
$\ps(\pi(w_j)) = \ps(\pi(w_{-j}))$. Further $\ps(e^{w_j \lambda -
  \lambda}) = v^{(\lambda-w_j \lambda,\rho)} = v^{(\lambda, \rho -
  w_j^{-1}\rho)}$.  Observe also that $w_j^{-1}$ is either $w_j$ or
$w_{-j}$ depending on whether $j$ is odd or even.
The above remarks together with equation \eqref{dlpl} imply that
$$\ps(W_\lambda(t) e^{-\lambda} \dtl \hl) = \sum_{j \in \integers}
v^{(\lambda, \rho - w_j \rho)} \ps(\pi(w_j)).$$

\onehalfspace
Let $\lambda$ be of {\em level} $l \geq 0$, i.e, $(\lambda, \delta)=l$. If
$(\lambda, \alpha_0)=p$, one has $0 \leq p \leq l$ and 
$(\lambda, \rho
- w_j \rho) = (\lambda, j \alpha_0 + {j \choose 2} \delta) = p j
 + l {j \choose 2} $.  Let $\fl[l]{p}:= \ps(W_\lambda(t) e^{-\lambda} \dtl \hl)$. 

\singlespace
Using the principal specialization
of $\pi(w_j)$, we obtain
\begin{equation} \label{maineqn}
 \fl[l]{p} = \sum_{j \in \integers} v^{ p
  j + l {j \choose 2}} \, t^j \, \frac{\hyp[j]{t^{\s -1} v ; v^2}}{\hyp[j]{t
      v ; v^2}} =   \sum_{j \in \integers} \left[ \frac{1+v^{(l-2p)j}}{2} \right] \, v^{ p
  j + l {j \choose 2}} \, t^j \, \frac{\hyp[j]{t^{\s -1} v ; v^2}}{\hyp[j]{t
      v ; v^2}} 
\end{equation}
where for the last equality, we used the fact that $\ps(\pi(w_j)) = \ps(\pi(w_{-j}))$.

If a closed form expression for $\fl[l]{p}$ can be found, equations
\eqref{psdtl} and \eqref{maineqn} allow us to 
determine the principal specialization of
$\hl$ via the relation 
\beq \label{reln}
\ps(\ea[\lambda] \hl) = \dfrac{\hyp{tv, tv,
    tv^2; v^2}}{\hyp{v, v, v^2; v^2}}
\dfrac{\fl[l]{p}}{W_\lambda(t)}.
\eeq
We now show that the sum defining $\fl[l]{p}$ can be written as an
explicit infinite product for certain special values of $(l,p)$.
These are $$(l,p) = (0,0), (1,0), (1,1), (2,0), (2,1), (2,2), (4,1),
(4,3).$$  We use classical summation theorems for bilateral
$q$-hypergeometric series for this purpose (see Appendix).  Note that the symmetry
between $\alpha_0$ and $\alpha_1$ ensures that $\fl[l]{\,p} =
\fl[l]{\,l-p}$; it is thus enough to consider $0 \leq p \leq [l/2]$.
When $(l,p)$ equals $(0,0)$ or $(1,0)$, we recover the identities of
Macdonald \cite{igmgod} and Fishel-Grojnowski-Teleman
\cite{fgt,groj} respectively.

\subsection{} 
$l=0, p=0$: In this case, as shown by Macdonald \cite{igmgod}, $\fl[0]{0}$
is just the $\bil[1]$ bilateral $q$-hypergeometric series. 
Taking $a=t^{\s -1}v, b=tv, w=v^2, z=t$ in equation \eqref{appendixeqn1} of
the Appendix, and using equation \eqref{reln}, one
obtains (cf \cite[equation (4.3)]{igmgod}):
\begin{equation} \label{l0}
\ps(\ea[\lambda] \hl) = \frac{\hyp{t^2 v^2; v^2}}{\hyp{t v^2; v^2}}                          
\end{equation}

\smallskip
\subsection{}
$l=1, p=0$:
$$\fl[1]{0} = \sum_{j \in \integers} v^{j \choose 2} \, \frac{t^j
  \hyp[j]{t^{\s -1} v ; v^2}}{\hyp[j]{t v ; v^2}} \, \left[  \frac{1+v^j}{2} \right] = \sum_{j \in
      \integers} v^{j \choose 2} \, t^j \, \frac{\hyp[j]{\sqrt{t^{\s -1} v},
        -\sqrt{t^{\s -1} v} ; v}}{\hyp[j]{\sqrt{t v}, -\sqrt{t v} ;
          v}} \left[  \frac{1+v^j}{2} \right].$$  This can be summed by specializing Bailey's
        $\bil[6]$ identity (equation \eqref{appendixeqn6}) at $w=v,
        b = \sqrt{t^{\s -1} v}, c = - \sqrt{t^{\s -1} v},
  d=a^{\frac{1}{2}}$ and letting $e \to
        \infty, a \to 1$. The resulting identity is
\begin{equation} \label{l1}
\ps(\ea[\lambda] \hl) =  \frac{\hyp{t^2 v^2; v^2}}{\hyp{v; v^2}}    
\end{equation}
The unspecialized version of this identity is implicit in \cite{fgt,groj,svis-kosfoulkes}.

\smallskip
\subsection{}
$l=2, p=1$: 
In this case, we first set $u = v^2$.
$$\fl[2]{1} = \sum_{j \in \integers} u^{j \choose 2} \, (tu^{\s 1/2})^j \,
\frac{\hyp[j]{t^{\s -1} u^{\s 1/2} ; u}}{\hyp[j]{t u^{\s 1/2} ; u}}.$$ We
    recognize this as the specialization of the $\bil[6]$
    series at $w=u,  b = t^{\s -1}
    u^{\s 1/2}, c=a^{\frac{1}{2}}, d=-c$ and $e \to \infty, a
    \to 1$. We thereby obtain
\begin{equation} \label{l2}
\ps(\ea[\lambda] \hl) = \frac{\hyp{t v^2; v^2} \, \hyp{t^2 v^2;
    v^4}}{\hyp{v, v, -v^2; v^2}}. 
\end{equation}

\smallskip
\subsection{}
$l=2, p=0$: 
Again, with $u = v^2$, we have

$$\fl[2]{0} =  \sum_{j \in \integers} u^{j \choose 2} \, \frac{1 + u^{j }}{2} \, \frac{t^j \, \hyp[j]{t^{\s -1} u^{\s 1/2} ; u}}{\hyp[j]{t u^{\s 1/2} ; u}}$$
We now apply Bailey's $\bil[6]$ with $w=u, b=t^{\s -1} u^{\s 1/2}, c=
- a^{\frac{1}{2}} \, u^{\s 1/2}, d=a^{\frac{1}{2}}$ and let $e \to
        \infty, a \to 1$ to get:
\begin{equation} \label{l2-1}
\ps(\ea[\lambda] \hl) =  \frac{\hyp{t v; v} \, \hyp{-tv^2;
    v^2}}{\hyp{v, v, -v; v^2}} 
\end{equation}

\smallskip
\subsection{}
$l=4, p=1$: 
Setting $u = v^2$ as above gives
$$\fl[4]{1} = \sum_{j \in \integers} u^{2{j \choose 2}} \, (tu^{\s
  1/2})^j \,
\frac{\hyp[j]{t^{\s -1} u^{\s 1/2} ; u}}{\hyp[j]{t u^{\s 1/2} ; u}} \left[  \frac{1+u^j}{2} \right].$$ This
    is again a specialization of the $\bil[6]$ identity, this time at $w=u, b = t^{\s -1}
      u^{\s 1/2}, c = a^{\frac{1}{2}}$ and $e \to \infty, a \to 1$.  Thus:
\begin{equation} \label{l4}
\ps(\ea[\lambda] \hl) = \frac{\hyp{t v; v}}{\hyp{v, v; v^2}}  
\end{equation}

\noindent
{\bf Remark 1:} (a). Besides the already known cases of levels 0 and 1, we
have thus obtained infinite product expressions for $\ps(\ea[\lambda]
\hl)$ for all dominant weights of level 2 and some of level 4. A 
natural question is whether these exhaust the cases where such infinite
product expressions exist. It is easy to see from equation \eqref{maineqn} that
for any dominant weight $\lambda$, 
further specializing $t=v^g$ for an odd natural number $g$, reduces 
$\fl[l]{p}$ to a rational function of $v$. When $\fl[l]{p}$ admits
an infinite product expression involving terms of the 
form $\hyp{\pm t^i
  \, v^j; \, v^k}^{{}^{}}$ as above, all zeros and
poles of these rational functions have modulus 1. 
Computational data (using {\small MAPLE}) and this
latter observation suggest that the only
 $\lambda$'s of level $\leq 10$ which admit
such infinite product expressions are the ones we have already found.  

\vspace{2mm}
\noindent
(b). We also recall that at $t=0$ and $t=1$, $\ea[\lambda] \hl$ reduces to
$\ea[\lambda] \sch$ and $\ea[\lambda] m_{\s \lambda}$ respectively. 
The principal specializations of both these have 
well known infinite product expressions {\em for all} $\lambda
\in P^+$ (by \cite[(10.10.1)]{kac} and the Jacobi triple product
identity respectively). Given $\lambda \in P^+$ with $(\lambda,
\delta) = l > 0$, $(\lambda, \alpha_0) =p$, one has
\begin{equation} \label{schps}
\ps(e^{-\lambda} \sch) = \frac{\hyp{v^{p+1}, v^{l-p+1}, v^{l+2};
    \, v^{l+2}}}{\hyp{v, v, v^2; \, v^2}} 
\eeq
\beq \label{orbsumps}
\ps(e^{-\lambda} m_{\s \lambda}) = \frac{1}{\# W_\lambda} \hyp{-v^p, -v^{l-p},
  v^l; \, v^l}
\eeq
For our special $\lambda$'s, equations \eqref{l1}-\eqref{l4} interpolate
between these two expressions (the elementary identity $\hyp{-w; w} = 1/ \hyp{w; w^2}$ is useful in checking
 this explicitly).

\section{Main theorems}\label{mainthms}
We now use the principal specializations of \S \ref{prspsec} to
determine the $t$-string function $\alm$ for the special $\lambda$'s of
level 2 and 4 and all $\mu \in \Max(\lambda)$. 
For ease of notation, let $\Lambda_0, \Lambda_1$ denote (any choice
of) fundamental weights, i.e, satisfying $(\Lambda_i, \alpha_j) =
\delta_{ij}$ for $i, j=0, 1$.
We recall that explicit expressions for 
the level 0 and 1 \tsfs $\aoo$ and $\alolo$ are known \cite{igmgod,fgt,groj,svis-kosfoulkes}: 
$$ a^{\s 0}_{\s 0}(t,v) =  \dfrac{\hyp{tq;q}}{\hyp{t^2 q; q}}    \;\;\;\;\;\;\;\;\; \text{ and } \;\;\;\;\;\;\;\;\; a^{\s \Lambda_0}_{\s \Lambda_0}(t,v) =  \dfrac{1}{\hyp{t^2 q; q}}.$$
These are respectively equivalent to the Macdonald and Macdonald-Mehta constant term
identities for $\afsl$. In Theorem \ref{chthm} below, we state 
these constant term identities. But first some notation:
let $$\chke:=\prod_{\alpha \in \rrep}
\frac{1-e^{-\alpha}}{1-te^{-\alpha}} = \frac{\hyp{\ea[\alpha_1],\,
  q \, e^{\alpha_1}; \, q}}{\hyp{t \ea[\alpha_1],\,
  t q \, e^{\alpha_1}; \, q}} \text{ (the Cherednik kernel)} $$ and let  
$\Thetaone(e^{-\alpha_0}, e^{-\alpha_1}) = \sum_{j \in \integers}
(e^{\s -\alpha_0})^{\s j^2}  (e^{\s -\alpha_1})^{\s j^2 -j}  =\sum_{j \in
  \integers} e^{\s j \alpha_1}  q^{j^2}$. We recognize $\Thetaone =
\Thetaone(e^{-\alpha_0}, e^{-\alpha_1})$ as the theta function of the root lattice of the underlying finite
dimensional simple Lie algebra $A_1 \cong sl_2$.

\begin{theorem}\label{chthm}
{\em (Cherednik)} For the affine Lie algebra $\afsl$, one has \\
\be
\item $\ct(\chke) = \dfrac{\hyp{tq;q}^2}{\hyp{t^2 q; q} \hyp{q; q}}$

\vspace{2mm}
\item $\ct(\chke \, \Thetaone) = \dfrac{\hyp{tq;q}}{\hyp{t^2 q; q}}$ \qed
\ee
\end{theorem}

\noindent
For later use, we also define $\dtli := \prod_{k \geq 1}
\left(\frac{1-e^{-k\delta}}{1-te^{-k\delta}}\right) = \frac{\hyp{q;
    q}}{\hyp{tq; q}}$. Observe that  
\begin{equation} \label{mueqn}
\dtl = \chke \, \dtli \text{   and   }
\ct(\dtli) = \dtli.
\end{equation}
In the next three subsections, we derive closed-form expressions for
our \tsfs of levels 2 and 4.  As mentioned in the introduction, these
in turn determine the corresponding affine Hall-Littlewood functions
and give us new constant term identities of Cherednik-Macdonald-Mehta
type.

\subsection{{\bf $\lambda = \Lambda_0 + \Lambda_1$}} 
Then $\lambda$ is of level 2 and $\Max(\lambda) = \{\lambda\}$. Thus 
\beq \label{schl}
\sch =  a^{\s \lambda}_{\s \lambda}(t,q) \hl.
\eeq
 Applying the principal
  specialization, we get $\ps(\ea[\lambda] \, \sch) =
  a^{\s \lambda}_{\s \lambda}(t,v^2) \, \ps(\ea[\lambda] \, \hl)$. By equation
  \eqref{schps}, 
$\ps(\ea[\lambda] \, \sch) =  \frac{\hyp{v^2; \, v^4}}{\hyp{v, v ;\,  v^2}}$. 
Together with equation \eqref{l2}, this gives us the
  following.
\begin{theorem} \label{al0l1thm}
$$a^{\s \Lambda_0 + \Lambda_1}_{\s \Lambda_0 + \Lambda_1}(t,q) =
  \frac{1}{\hyp{tq;q}\hyp{t^2q;q^2}}$$ \qed
\end{theorem}

To obtain the corresponding constant term identity, first observe that
$a^{\s \lambda}_{\s \lambda} (t,q) = \ct(e^{-\lambda}\, \dtl\, \sch) =
a^{\s \lambda}_{\s \lambda} (1,q)  \dtli \, \ct(\chke \, e^{-\lambda}
m_{\s \lambda})$. 
The first equality is by equation \eqref{almeq}, and the second
follows from equation \eqref{mueqn} and equation \eqref{schl} at $t=1$.
For $\lambda = \Lambda_0 + \Lambda_1$, it is straightforward to see that
$$e^{-\lambda} m_{\s \lambda} =\sum_{\mu \in W.\lambda}
e^{\mu - \lambda} = \sum_{j \in \integers}
(e^{-\alpha_0})^{\frac{j(j+1)}{2}} \,
    (e^{-\alpha_1})^{\frac{j(j-1)}{2}} = \Thetaram(\ea[\alpha_0],
  \ea[\alpha_1])$$ where $\Thetaram$ is the Ramanujan Theta function
  \cite[p. 34]{berndt}. Observe that $\Thetaram$ is also equal to
  $\sum_{j \in \integers} q^{j \choose 2} (e^{-\alpha_1}) ^j$.
Theorem \ref{al0l1thm} thus implies the
following level 2 analog of Cherednik's difference
Macdonald-Mehta identity, for $\afsl$.
\begin{corollary} \label{ct-lev2-1}
$$\ct(\chke \, \Thetaram) = \frac{\hyp{q;q^2}}{\hyp{t^2 q; q^2}}$$ \qed
\end{corollary}

\noindent
Next, recalling the well known fact that $\sch[\s \Lambda_0 + \Lambda_1] = e^{\s \Lambda_0 + \Lambda_1} 
\, \prod_{\alpha \in \rp} (1 + e^{\s -\alpha})$ for \afsl \cite[Ex 10.1]{kac}, 
 one observes that equation \eqref{schl} and theorem \ref{al0l1thm} together determine the affine Hall-Littlewood 
function $\hl[\Lambda_0 + \Lambda_1]$:
\begin{corollary}
$$\hl[\Lambda_0 + \Lambda_1] = e^{\Lambda_0 + \Lambda_1} \, \hyp{-e^{-\alpha_1}, -qe^{\alpha_1}, -q, tq; \, q} \hyp{t^2q;q^2} $$
\qed
\end{corollary}

\smallskip
\subsection{{\bf $\lambda = 2 \Lambda_0$}} Again, $\lambda$ is  
of level 2, with $\Max(\lambda) = \{2\Lambda_0, 2\Lambda_0 -
\alpha_0\}$. Thus 
\begin{equation} \label{2locase}
\sch[2\Lambda_0] =  a^{\s 2\Lambda_0}_{\s 2\Lambda_0}(t,q) \, \hl[2\Lambda_0] +
  a^{\s 2\Lambda_0}_{\s 2\Lambda_0 - \alpha_0}(t,q) \, \hl[2\Lambda_0 -
    \alpha_0]
\end{equation}
\onehalfspace
 We observe that $(2\Lambda_0 - \alpha_0, \alpha_0) =0
  $ and $(2\Lambda_0 - \alpha_0, \alpha_1) =2
 $. The $\alpha_0 \leftrightarrow \alpha_1$ symmetry implies
  $\ps(e^{-2\Lambda_0} \hl[2\Lambda_0]) =  \ps(e^{-(2\Lambda_0 - \alpha_0)}
  \hl[2\Lambda_0 - \alpha_0])$. Principally specializing equation
  \eqref{2locase}, one obtains
\begin{equation} \label{2lops}
\ps(\ea[2\Lambda_0] \, \sch[2\Lambda_0]) = \left(
a^{\s 2\Lambda_0}_{\s 2\Lambda_0}(t,v^2)  + v\, a^{\s 2\Lambda_0}_{\s 2\Lambda_0 - \alpha_0}(t,v^2)
 \right)  \ps(\ea[2\Lambda_0] \hl[2\Lambda_0]).
\end{equation}
\singlespace
We now set $v \mapsto -v$, and use equations \eqref{l2-1}, 
\eqref{schps} to determine the individual \tsfs. 
The result is summarized in the following theorem. 

For ease of
notation, we let $\nh[a] := \hyp{a; q}$ in the rest of the article.
\begin{theorem} \label{a2l0thm}
\begin{enumerate}
\item  $a^{\s 2\Lambda_0}_{\s 2\Lambda_0}(t,q) = \dfrac{\nh[-tq^{\frac{1}{2}}] +
  \nh[tq^{\frac{1}{2}}]}{2\nh[t^2 q]}$
\item $a^{\s 2\Lambda_0}_{\s 2\Lambda_0 - \alpha_0}(t,q) =
  \dfrac{\nh[-tq^{\frac{1}{2}}] -
  \nh[tq^{\frac{1}{2}}]}{2q^{\frac{1}{2}} \, \nh[t^2 q]}$ \ee
  \qed
\end{theorem}

\noindent
This theorem allows us to derive two further constant term
identities. To state these, we introduce $\Thetatwo =
\Thetatwo(e^{-\alpha_0}, e^{-\alpha_1}) := \sum_{j \in \integers}
(e^{-\alpha_0})^{2 j^2} \,    (e^{-\alpha_1})^{2j^2 - 2j}$. Observing
that $\Thetatwo = e^{-2\Lambda_0} m_{\s 2 \Lambda_0}$, it is clear that 
$\Thetatwo(e^{-\alpha_0}, e^{-\alpha_1}) = \Thetaone (e^{-2\alpha_0},
e^{-2\alpha_1})$. We also let $\Thetatwobar :=
\Thetatwo(e^{-\alpha_1}, e^{-\alpha_0})$.
With this notation, we claim the following.
\begin{corollary} \label{ctlev2}
\begin{align}
{\mathrm (1)} \;\;\;\;\; \;\;\;\;\;\;\; \ct(\chke \, \Thetatwo) &=
  \frac{q}{1-q}\;\frac{\nh[tq]}{\nh[t^2q]}\left(
  \frac{\nh[-tq^{\frac{1}{2}}]^{{}^{\mathrm{odd}}} \, q^{-\frac{1}{2}}
  - \nh[-tq^{\frac{1}{2}}]^{{}^{\mathrm{even}}}}
  {\nh[-q^{\frac{1}{2}}]^{{}^{\mathrm{even}}}}\right) \notag\\ &
  \notag\\ {\mathrm (2)} \;\;\;\;\, \ct(\chke \, \Thetatwo
  e^{-\alpha_1}) &= \frac{1}{1-q}\,\frac{\nh[tq]}{\nh[t^2q]}\left(
  \frac{\nh[-tq^{\frac{1}{2}}]^{{}^{\mathrm{even}}} - q^{\frac{1}{2}}
  \nh[-tq^{\frac{1}{2}}]^{{}^{\mathrm{odd}}}}{\nh[-q^{\frac{1}{2}}]^{{}^{\mathrm{odd}}}
  q^{-\frac{1}{2}}} \right)\notag
\end{align} 

\smallskip
\noindent
where for a series $\xi = \sum_{k=0}^{\infty} a_k \, q^{k/2}$, we let
$\xi^{{}^{\mathrm{odd}}} := \displaystyle\sum_{k \text{ odd}} a_k \, q^{k/2}$ \\and
$\xi^{{}^{\mathrm{even}}} := \displaystyle\sum_{k \text{ even}} a_k \, q^{k/2}$. 
\end{corollary}

\noindent
{\em Proof:} Setting $t=1$ in equation \eqref{2locase} gives:
\begin{equation} \label{2locase-1}
\sch[2\Lambda_0] =  a^{\s 2\Lambda_0}_{\s 2\Lambda_0}(1,q) m_{\s 2\Lambda_0} +
  a^{\s 2\Lambda_0}_{\s 2\Lambda_0 - \alpha_0}(1,q) m_{\s 2\Lambda_0 -
    \alpha_0}.
\end{equation}
Again, since $(2\Lambda_0 - \alpha_0, \alpha_0) =0
  $ and $(2\Lambda_0 - \alpha_0, \alpha_1) =2$, the $\alpha_0
\leftrightarrow \alpha_1$ symmetry argument
used before shows  $e^{-(2\Lambda_0 -\alpha_0)} m_{\s 2\Lambda_0 -
  \alpha_0} = \Thetatwobar$.  Equations \eqref{2locase-1}
and \eqref{almeq} thus imply the following relations.
\begin{align}
a^{\s 2\Lambda_0}_{\s 2\Lambda_0}(t,q) &= a^{\s 2\Lambda_0}_{\s 2\Lambda_0}(1,q) \ct(\dtl \Thetatwo) +
  a^{\s 2\Lambda_0}_{\s 2\Lambda_0 - \alpha_0}(1,q)  \ct(\dtl \ea[\s \alpha_0]
  \Thetatwobar) \label{one}\\
a^{\s 2\Lambda_0}_{\s 2\Lambda_0 -\alpha_0}(t,q) &= a^{\s
  2\Lambda_0}_{\s 2\Lambda_0}(1,q) \, q^{-1}\, \ct(\dtl e^{-\alpha_1} \Thetatwo) +
  a^{\s 2\Lambda_0}_{\s 2\Lambda_0 - \alpha_0}(1,q)  \ct(\dtl  \Thetatwobar)\label{two}
\end{align}
For $\beta = c_0 \alpha_0 + c_1 \alpha_1 \in Q$, let $\widehat{\beta} :=
c_0 \alpha_1 + c_1 \alpha_0$. Given $\xi = \sum_{\beta \in Q^+}
f_{\beta} \, e^{-\beta}$, define $\widehat{\xi} := \sum_{\beta \in Q^+}
f_{\beta} \, e^{-\widehat{\beta}}$. We observe the following easy
properties of the hat operation:
 (i) $\ct(\xi) = \ct(\widehat{\xi}\,)$, (ii) $\widehat{\xi \eta} = \widehat{\xi}\, \widehat{\eta}$ and (iii)
$\dtl$ is invariant under the hat operation.  These imply that $\ct(\dtl \Thetatwo)  =
\ct(\dtl  \Thetatwobar)$ and $ \ct(\dtl \, \ea[\alpha_0]
  \Thetatwobar) = \ct(\dtl \, e^{-\alpha_1} \Thetatwo)$. So, equations
  \eqref{one} and \eqref{two} above form a $2 \times 2$ system :
$$\begin{bmatrix}  \zeta_1 \\ \zeta_2
  \end{bmatrix} =
  {\textstyle A} \begin{bmatrix}\scriptstyle \ct(\dtl \Thetatwo) \\  \scriptstyle \ct(\dtl
    \, e^{-\alpha_1} \Thetatwo) \end{bmatrix}$$ where ${\zeta}_1 = 
 a^{\s 2\Lambda_0}_{\s   2\Lambda_0}(t,q), \; {\zeta}_2 = a^{\s
   2\Lambda_0}_{\s 2\Lambda_0 -\alpha_0}(t,q)$ and 
  the matrix 
  ${\textstyle A}$ is given by
$$A = \begin{bmatrix} a^{\s 2\Lambda_0}_{\s 2\Lambda_0}(1,q) &
    a^{\s 2\Lambda_0}_{\s 2\Lambda_0 - \alpha_0}(1,q) \\ 
 & \\ q^{\s -1} a^{\s
      2\Lambda_0}_{\s 2\Lambda_0}(1,q) &
    a^{\s 2\Lambda_0}_{\s 2\Lambda_0 - \alpha_0}(1,q) \end{bmatrix}.$$
Inverting $A$, and using equation \eqref{mueqn} completes the proof
of corollary \ref{ctlev2}. \qed

As before, theorem \ref{a2l0thm} can be used to derive an explicit
expression for $\hl[2\Lambda_0]$. One now uses equation \eqref{2locase}, and the
$\alpha_0 \leftrightarrow \alpha_1$ symmetry. We content ourselves
with stating the result of this calculation.

\begin{corollary}
With notation as above, one has:
$$ \hl[2\Lambda_0] = e^{2\Lambda_0} \, (\alpha \,\Thetatwo + \beta e^{-\alpha_0}\,
\Thetatwobar)$$ with $\alpha = \frac{1}{2} \left(
\frac{\nh[-q^\frac{1}{2}]}{\nh[-t q^\frac{1}{2}]} +
\frac{\nh[q^\frac{1}{2}]}{\nh[t q^\frac{1}{2}]} \right)  \,
\frac{\nh[t^2q]}{\nh[q]^{}}$ and $\beta = \frac{1}{2q^\frac{1}{2}} \left(
\frac{\nh[-q^\frac{1}{2}]}{\nh[-t q^\frac{1}{2}]} -
\frac{\nh[q^\frac{1}{2}]}{\nh[t q^\frac{1}{2}]} \right) \,
\frac{\nh[t^2q]}{\nh[q]^{}}$.
\qed
\end{corollary}

\medskip 
\subsection{{\bf $\lambda = 3 \Lambda_0 + \Lambda_1$}} In this case, $\lambda$ is  
of level 4, with $\Max(\lambda) = \{3\Lambda_0 + \Lambda_1, 3\Lambda_0 + \Lambda_1 -
\alpha_0\}$. We let $\Thetafour := e^{-(3 \Lambda_0 + \Lambda_1)}\,
m_{\s 3 \Lambda_0 + \Lambda_1}$  Reasoning as in the previous subsection and using the
principal specialization from equation \eqref{l4}, one obtains the
following theorem and its corollary.

\begin{theorem}\label{alev4thm}
\be
\item
$a^{\s 3 \Lambda_0 + \Lambda_1}_{ \s 3 \Lambda_0 +  \Lambda_1}(t,q)
  = \nh[-tq] \; a^{\s 2\Lambda_0}_{\s 2\Lambda_0}(t,q)$.

\vspace{1.5mm}
\item 
$a^{ \s 3 \Lambda_0 + \Lambda_1}_{ \s 3 \Lambda_0 +  \Lambda_1 - \alpha_0}(t,q) = \nh[-tq] \; a^{\s 2\Lambda_0}_{\s 2\Lambda_0-\alpha_0}(t,q)$.   \qed
\ee
\end{theorem}

\noindent

\begin{corollary} \label{ct-lev4}
\be
\item $\ct(\chke \, \Thetafour) = \dfrac{\nh[-tq]}{\nh[-q]}
  \ct(\chke \, \Thetatwo)$.
\item 
$\ct(\chke \, \Thetafour \, e^{-\alpha_1}) = \dfrac{\nh[-tq]}{\nh[-q]}
  \ct(\chke \, \Thetatwo \, e^{-\alpha_1})$.
\ee \qed
\end{corollary}

\noindent
As in level 2, there is an explicit formula for $\hl[3\Lambda_0 +
  \Lambda_1]$ as well; the details are omitted.

\vspace{1.5mm}
\noindent
{\bf Remark 2:} We observe that the $\bil[6]$ identity was central to
all our $\afsl$ \tsf computations. A generalization of our approach to the higher
rank affines $\kma = A_n^{(1)}$ might be possible by using a suitable
multivariable generalization of the $\bil[6]$ sum.
Several such (distinct) generalizations are known \cite{gustafson1,
  gustafson2, schlosser1, schlosser2, 
  vandiejen}, and it is an interesting question whether one of these
choices leads to explicit formulas for affine Hall-Littlewood
functions of small level, for $A_n^{(1)}$.

\section{Appendix}
For quicker reference and to fix notation, we give below the classical summation
formulas for bilateral basic hypergeometric series that have been used
in this article. 

For parameters $a_i, b_j$, recall that the series ${}_m \psi_n$ is
defined by
\beq
{}_m \psi_n \left( \begin{smallmatrix} a_1 \; a_2 \; \cdots\; a_m
  \\ b_1 \; b_2 \; \cdots \; b_n \end{smallmatrix} ; \; w, z \right)
:= \sum_{j \in \integers} \frac{\hyp[j]{a_1, a_2, \cdots, a_m; \,
    w}}{\hyp[j]{b_1, b_2, \cdots, b_n; \, w}} \, z^j
\eeq 
Ramanujan's $\bil[1]$ and Bailey's $\bil[6]$ sums  \cite[Appendix II]{gr} are given
below.
\begin{align}
\bil[1] \left( \begin{smallmatrix} a \\ b \end{smallmatrix} ; \; w, z
\right) &= \frac{\hyp{w, \frac{b}{a}, az, \frac{w}{az}; \, w}}{\hyp{b,
    \frac{w}{a},  z, \frac{b}{az}; \, w}}\label{appendixeqn1}\\
\bil[6] \left(\begin{smallmatrix} wa^{\frac{1}{2}} & -wa^{\frac{1}{2}} & b & c &
  d & e\\ a^{\frac{1}{2}} & -a^{\frac{1}{2}} & \frac{aw}{b} &
  \frac{aw}{c} & \frac{aw}{d} & \frac{aw}{e} \end{smallmatrix} ; \; \begin{matrix} w,
  \, \frac{wa^{2}}{bcde} \end{matrix} \right)
  &= \frac{\hyp{aw, \frac{aw}{bc}, \frac{aw}{bd}, \frac{aw}{cd},
      \frac{aw}{be}, \frac{aw}{ce}, \frac{aw}{de}, w, \frac{w}{a}\, ; \,
      w}}{\hyp{\frac{aw}{b}, \frac{aw}{c}, \frac{aw}{d}, \frac{aw}{e},
      \frac{w}{b}, \frac{w}{c}, \frac{w}{d}, \frac{w}{e}, 
       \frac{wa^2}{bcde} \, ; \, w}}\label{appendixeqn6}
\end{align}

\end{document}